%% file: word.tex
\documentclass[12pt]{amsart}

\input{header.tex}

\newcommand{\produces}{\Implies} 

\newcommand{\height}{{\operatorname{height}}} 
\newcommand{\rev}{{\operatorname{rev}}} 

\newcommand{\Split}{{\operatorname{Split}}} 
\newcommand{\Compact}{{\operatorname{Compact}}} 
\newcommand{\SimpleCompact}{{\operatorname{SimpleCompact}}} 

\newcommand{\ova}{{\overline{a}}}  
\newcommand{\ovb}{{\overline{b}}}  
\newcommand{\ovc}{{\overline{c}}}  
\newcommand{\ovd}{{\overline{d}}}  
\newcommand{\ove}{{\overline{e}}}  
\newcommand{\ovt}{{\overline{t}}}  

\begin{document}

\title{Polynomial-time word problems}

\author{Saul Schleimer}
\address{\hskip-\parindent
        Department of Mathematics\\
	Rutgers University\\
        Piscataway, New Jersey 08854}
\email{saulsch@math.rutgers.edu}
\urladdr{http://www.math.rutgers.edu/$\sim$saulsch}

\thanks{This work is in the public domain.}

\date{\today}

\begin{abstract}
We find polynomial-time solutions to the word problem for
free-by-cyclic groups, the word problem for automorphism groups of
free groups, and the membership problem for the handlebody subgroup of
the mapping class group.  All of these results follow from observing
that automorphisms of the free group strongly resemble straight line
programs, which are widely studied in the theory of compressed data
structures.  In an effort to be self-contained we give a detailed
exposition of the necessary results from computer science.
\end{abstract}

\maketitle

\section{Introduction}

Automorphisms of the free group are closely connected to two
techniques in computer science: {\em string matching} and {\em
compression}.  The relevance of the first is obvious.  The second is
less clear.  So, consider the fact that an automorphism of complexity
$n$ can produce, by acting on a generator, words of size at most
$\exp(n)$.  Now, there are only $\exp(n)$ such automorphisms while
there are $\exp(\exp(n))$ words available as output.  Thus most words
in the free group {\em cannot} be obtained in this way.  Those which
can are highly regular and thus susceptible to compression.

Compression techniques have already made an appearance in algorithmic
topology.  Word equations play a starring role in the work of
Schaefer, Sedgwick, and Stefankovic~\cite{SchaeferEtAl02}.  One of the
problems they consider, connectedness of normal curves and surfaces,
is also addressed by the orbit-counting techniques of Agol, Hass, and
Thurston~\cite{AgolEtAl02}.  Both results rely, directly or
indirectly, on Plandowski's Algorithm~\cite{Plandowski94}
(Theorem~\ref{Thm:PlandowskiI} below).

The structure of the paper is as follows: Section~\ref{Sec:SLP}
reviews {\em straight line programs} and also a slight variant, {\em
composition systems}.  Such programs are called {\em compressed
words}.  Section~\ref{Sec:Lohrey} is an exposition of Lohrey's
Theorem~\cite{Lohrey04}:

\begin{restate}{Theorem}{Thm:WordProbInFree}
The word problem for compressed words in the free group is solvable in
polynomial time.
\end{restate}

The same techniques solve several open problems: we give
polynomial-time algorithms for the word problem in free-by-cyclic
groups (Section~\ref{Sec:FreeByCyclic}), the word problem in
$\Aut(F_m)$ (Section~\ref{Sec:Aut}), and various membership problems
(Section~\ref{Sec:Membership}).  In each case the compression
technique ``accelerates'' the obvious exponential-time algorithm.

Sections~\ref{Sec:Hagenah} and~\ref{Sec:Plandowski} 
give expositions of theorems due to Hagenah~\cite{Hagenah00} and
Plandowski~\cite{Plandowski94}, upon which Lohrey's Theorem relies.
As a result the computer science portions of the paper are
self-contained.  I have been somewhat more permissive when using
techniques from combinatorial group theory or three-manifolds.

The paper ends with a brief appendix (Section~\ref{Sec:Surfaces})
indicating how these techniques extend to closed surface groups.

\subsection*{Acknowledgements} 

I thank Martin Bridson and Nathan Broaddus for both inspiring seminar
talks and for illuminating conversations.  Without such models this
paper would not exist.  I also thank Ilya Kapovich for his comments on
an early version of this paper.

\section{Straight line programs and composition systems}
\label{Sec:SLP}

Recall that if $\calL$ is a set of characters, or an {\em alphabet},
then the {\em Kleene closure} $\calL^*$ is the set of all {\em words}
(finite strings of characters) in $\calL$.  If $w \in \calL^*$ then
$|w|$ is the {\em length} of $w$: the number of characters of $w$
counted with multiplicity.  We write $\epsilon$ for the {\em empty
string}: the word of length zero.  If $u$ and $v$ are words we write
$u = v$ if and only if $u$ and $v$ are identical as strings.

Fix a word $w$ and assume $0 \leq i \leq j \leq |w|$.  We take
$w[i:j]$ to be the substring starting immediately before the $i+1^\th$
character of $w$ and ending just after the $j^\th$ character.  Thus
$w[i:i] = \epsilon$ and, in general, $|w[i:j]| = j - i$.  By
convention negative indices count from the end of $w$.  Thus:
$$w[-j:-i] = w[|w|-j:|w|-i].$$ 
The following abbreviations will be useful: $w[:i]$ for $w[0:i]$,
$w[i:]$ for $w[i:|w|]$, and $w[i]$ for $w[i:i+1]$.  Thus $w[0]$ is the
first character of $w$ while $w[-1]$ is the last.

If $u, v$ are words in $\calL^*$ then their concatenation is denoted
by $u \cdot v$.  This leads to the pleasant identity $w = w[:i] \cdot
w[i:]$.  Also, in order to {\em rotate} a word $w$ exactly $i$
characters simply form $w[i:] \cdot w[:i]$.  For any word $w \in
\calL^*$ define $\rev(w)$ to be the {\em reverse} of $w$: so
$\rev(w)[i] = w[-i - 1]$.

A {\em straight line program} $\AA = \langle \calL, \calA, A_n, \calP
\rangle$ contains the following: a finite alphabet $\calL = \{ a_1,
\ldots, a_m \}$ of {\em terminal} characters, a disjoint finite
alphabet $\calA = \{ A_1, \ldots, A_n \}$ of {\em non-terminal}
characters, a {\em root} terminal $A_n \in \calA$, and a set $\calP =
\{ A_i \produces W_i \}$ of {\em production rules}.  These last allow
us to replace a non-terminal $A_i$ with its {\em production}: a
(possibly empty) word $W_i$ in $(\calL \cup \calA)^*$.  Every
non-terminal $A_j$ appearing in $W_i$ has index $j < i$.  When indices
are unimportant we shall use the non-terminal $A$ to represent the
root of $\AA$.

To justify the term ``root'' define the {\em production tree} of a
character in $\calL \cup \calA$ as follows: The tree for a terminal $a
\in \calL$ is a single vertex, labelled $a$.  The tree for a
non-terminal $A_i \in \calA$ is planar, has root labelled $A_i$, and,
attached to the root in left-to-right order, there is a copy of the
production tree for every character of $W_i$.

If $B \in \calL \cup \calA$ then the {\em height} of $B$, denoted
$||B||$, is the height of the production tree of $B$: the maximal
distance from the root to a leaf.  For example, terminals have height
zero.

Define $w(A) = w_A$ to be the word of $\calL^*$ that results from
running the straight line program $\AA$.  That is, if $A_n$ is the
root then produce $W_n$, replace all non-terminals appearing with
their productions, and continue doing so until the resulting word lies
in $\calL^*$.  This is exactly the word appearing at the leaves of the
production tree of $\AA$.  Also, define $w(A_i) = w_{A_i}$ to be the
word in $\calL^*$ produced by the non-terminal $A_i$.

A straight line program is in {\em Chomsky normal form} if every
production $W_i$ has length one or two: all of the former lie in
$\calL$ while all of the latter lie in $\calA^*$.

\begin{remark}
\label{Rem:NormalForm}
Every straight line program can be placed in normal form in time
polynomial in $\Sigma_i |W_i|$.  This is done by introducing dummy
non-terminals.
\end{remark}

\begin{remark}
\label{Rem:Time}
Two remarks are in order about the meaning of {\em polynomial time}.
First, the variable in question is the {\em bit-size} of the input.
For example, a straight line program in normal form, with $n$
non-terminals, has bit-size $O(n \log_2(n))$.  
The exact bit-size depends on a choice of {\em encoding}.  Second, the
precise running-time of an algorithm depends on the {\em model of
computation}.  However, it is generally true that changing encoding or
model of computation transforms the bit-size or running-time by a
polynomial function.  Thus the claim that some problem may be solved
in polynomial time is essentially independent of these choices.
\end{remark}

\begin{example}
\label{Exa:Fib}
Here is the canonical example of a straight line program:
$$\FF = 
\left\langle 
\begin{array}{l}
  \{ a, b \}, \{ F_i \}, F_n,  \\
  \{ F_{i} \produces F_{i-1} \cdot F_{i-2} \}_{i = 3}^{n} \cup 
  \{ F_2 \produces a, F_1 \produces b \} \\
\end{array}
\right\rangle.$$
\noindent
So $w(F_i)$ is the $i^\th$ Fibonacci word.  For example
\begin{eqnarray*}
w(F_8) & = & abaababaabaababaababa, \\
w(F_9) & = & abaababaabaababaababaabaababaabaab.
\end{eqnarray*}
It follows that the length $|w(F_n)|$ grows exponentially with $n$.   
\end{example}






{\em Composition systems} are a more flexible version of straight line
programs (in normal form), introduced in~\cite{GKPR96}: If $A, B, C
\in \calA$, then productions of the form $A \produces B[i:j] \cdot
C[k:l]$ are allowed.  Here $B[i:j]$ is a {\em truncated} non-terminal.
Truncated non-terminals only appear on the right hand side of
productions, never on the left.  A truncated non-terminal $B[i:j]$ is
well formed if the indices satisfy $0 \leq i \leq j \leq |w_B|$.
Define $w(B[i:j]) = w_B[i:j]$.
Repeated truncation behaves quite simply: $(B[i:j])[k:l] =
B[i+k:i+l]$.


\begin{example}
\label{Exa:Composition}
Consider the straight line program in normal form
$$\AA = 
\left\langle 
\begin{array}{l}
  \{a, b\}, \{ A_i \}, A_n,  \\
  \{ A_{i} \produces A_{i-1} \cdot A_{i-1} \}_{i = 4}^{n} \cup \{ A_3
  \produces A_2 \cdot A_1, A_2 \produces a, A_1 \produces b \} \\
\end{array}
\right\rangle.$$
Of course, $w(A_i) = (ab)^{2^{i-3}}$ for $i \geq 3$.  Truncate a
little to obtain:
$$\BB = 
\left\langle 
\begin{array}{l}
  \{a, b\}, \{ B_i \}, B_n,  \\
  \{ B_{i} \produces B_{i-1}[1:] \cdot B_{i-1}[1:] \}_{i = 5}^n \cup \\
  \cup \{ B_4 \produces B_3 \cdot B_3, B_3 \produces B_2 \cdot B_1,
  B_2 \produces a, B_1 \produces b \} \\ 
\end{array} 
\right\rangle.$$
Now the output appears to be more interesting: for example
$$w(B_8) = bababbabbbababbabbababbabbbababbab.$$
\end{example}


For future use, we record:

\begin{lemma}
\label{Lem:LengthBound}
If $\AA$ is a composition system of height $||A||$ then $|w_A| \leq
2^{||A||}$.  \qed
\end{lemma}


\begin{lemma}
\label{Lem:ComputeLength}
There is a polynomial-time algorithm that, given a composition system
$\AA$, computes $|w_B|$ for all $B \in \calA$. \qed
\end{lemma}

\begin{lemma}
\label{Lem:ComputeSubword}
There is a polynomial-time algorithm that, given a composition system
$\AA$ and an integer $i$, computes the character $w_A[i]$. \qed
\end{lemma}

\begin{remark}
\label{Rem:ComputeSubword}
In particular, suppose we are given $i, j$ with $0 \leq i \leq j \leq
|w_A|$.  Then we may compute the word $w_A[i:j]$ in time polynomial in
$j - i$ and in the size of $\AA$.
\end{remark}


It is perhaps surprising that the expressive power of composition
systems and straight line programs are nearly the same.  In Hagenah's
thesis~\cite[Chapter~8]{Hagenah00} we find: 

\begin{restate}{Theorem}{Thm:Hagenah}
There is a polynomial-time algorithm that, given a composition system
$\AA$, finds a straight line program $\XX$ with the property that $w_X
= w_A$.
\end{restate}

The algorithm and the proof of correctness are presented, for the
convenience of the reader, in Section~\ref{Sec:Hagenah}.  A more
subtle result is due to Plandowski~\cite{Plandowski94}:

\begin{restate}{Theorem}{Thm:PlandowskiI}
There is a polynomial-time algorithm that, given straight line
programs $\AA$ and $\XX$ in normal form, decides whether or not $w_A =
w_X$.
\end{restate}

In an attempt to be self-contained a proof appears in
Section~\ref{Sec:Plandowski}.  Gasieniec, Karpinski, Plandowski, and
Rytter~\cite{GKPR96} strengthen Theorem~\ref{Thm:PlandowskiI} as
follows:

\begin{theorem}
\label{Thm:PlandowskiII}
There is a polynomial-time algorithm that, given composition systems
$\AA$ and $\XX$, computes the largest integer $k \geq 0$ so that $w_A[:k]
= w_X[:k]$.
\end{theorem}

\begin{proof}
This follows from Theorems~\ref{Thm:Hagenah}
and~\ref{Thm:PlandowskiI}.  Compute $k$ via binary search.
\end{proof}


\section{Lohrey's algorithm}
\label{Sec:Lohrey}

We now turn our attention to the free group.  Let $\calL_m = \{ a_i,
\ova_i \}_{i = 1}^m$.  Let $\overline{\cdot} \from \calL_m \to
\calL_m$ be the obvious involution.  Given a word $w(a_i) \in
\calL_m^*$ define $\overline{w}$ to be $\rev(w(\ova_i))$.

{\em Compressed word} is the umbrella term for a straight line program
or composition system which produces a word in $\calL_m^*$.  The
involution given above extends to compressed words; constructions like
$\overline{B}$ are allowed on the right hand side of productions.  If
$A \produces B[i:j] \cdot C[k:l]$ then $\overline{A} \produces
\overline{C}[-l:-k] \cdot \overline{B}[-j:-i]$.  Finally, if $A
\produces a$ then $\overline{A} \produces \overline{a}$.  See
Example~\ref{Exa:Auto} below for an illustration.

\begin{lemma}
\label{Lem:ComputeInverse}
There is a polynomial-time algorithm that, given a compressed word
$\AA$, computes a new compressed word $\overline{\AA}$ so that
$w(\overline{A}) = \overline{w(A)}$. \qed
\end{lemma}

Now suppose that $\{ a_1, \ldots, a_m \}$ generate the free group
$F_m$ and that $\ova_i$ is the inverse of $a_i$.  Recall that a word
in the free group is {\em freely reduced} if it has no subwords of the
form $a_i\ova_i$ or $\ova_i a_i$ for any $i \in \{ 1, \ldots m \}$.  A
word is {\em cyclically reduced} if all of its rotations are freely
reduced.

\begin{example}
\label{Exa:Auto}
To simplify notation inside this example, let $\calL_2 = \{ a, b,
\ova, \ovb \}$.  Form a straight line program 
$$\AA = 
\left\langle 
\begin{array}{l}
\calL_2, \{ A_i, B_i \}_{i=0}^n, A_n, \\
  \{ A_{k+1} \produces B_k, B_{k+1} \produces B_k A_k \overline{B}_k
  \}_{k = 1}^{n - 1} \cup \{ A_0 \produces a, B_0 \produces b \} \\
\end{array}
\right\rangle.$$
\noindent
Thus
$$w(A_5) = 
ba\ovb b b\ova \ovb ba\ovb ba\ovb \ovb b\ova \ovb ba\ovb b b\ova
\ovb ba\ovb b b\ova \ovb b\ova \ovb ba\ovb \ovb b\ova \ovb.$$
Note the close relation with $\varphi^5(a)$ where $\varphi \from F_2
\to F_2$ is the automorphism $a \mapsto b, b \mapsto ba\ovb$.  Notice
also that $w(A_5)$ has both free and cyclic reductions.
\end{example}






Using the generalization of Plandowski's work (given as
Theorem~\ref{Thm:PlandowskiII} above), Lohrey~\cite{Lohrey04} has
proven:

\begin{theorem}
\label{Thm:Lohrey}
There is a polynomial-time algorithm that, given a straight line
program $\AA$ in normal form, finds a composition system $\XX$ with
the property that $w_X$ is the free reduction of $w_A$.
\end{theorem}



\begin{proof}
Induct on $n$.  Suppose $\AA = \langle \calL_m, \calA, A_n, \calP
\rangle$ is the given straight line program.  Now define a composition
system $\XX = \langle \calL_m, \calX, X_n, \calQ \rangle$.  For every
non-terminal $A_i$ of height one in $\calA$ place $X_i$ in $\calX$ and
add the production $X_i \produces w(A_i)$ to $\calQ$.

Now, if $\height(A_n) = 1$ then there is nothing more to prove.  So
assume that $\height(A_n) \geq 2$.  By induction assume that $X_i$,
for $i < n$, lies in $\calX$ and assume that the corresponding
production lies in $\calQ$.  Thus $w(X_i)$ is freely reduced for all
$i < n$.  Now place $X_n$ in $\calX$ and consider what $X_n$ will
produce.  

Suppose that $A_n \produces A_i \cdot A_j$.  Build $\overline{X}_i$
using the algorithm of Lemma~\ref{Lem:ComputeInverse}.  Apply the
algorithm of Theorem~\ref{Thm:PlandowskiII} to find the largest $k$ so
that $w(\overline{X}_i)[:k] = w(X_j)[:k]$.  Add the production
$$X_n \produces X_i[:-k] \cdot X_j[k:]$$ 
to $\calQ$.  The word $w(X_n)$ is now freely reduced.
\end{proof}


\begin{example}
\label{Exa:AutoCont}
We continue Example~\ref{Exa:Auto}.  Given $\AA$ as above Lohrey's
algorithm produces the following composition system:
$$\XX = 
\left\langle 
\begin{array}{l}
  \calL_2, \{ X_i, Y_i \}, X_n,  \\
  \{ X_{i+1} \produces Y_{i}, Y_{i+1} \produces Y_{i}[:2] \cdot
  \overline{Y}_{i} \}_{i = 1}^{n-1} \cup \\
  \cup \{ X_1 \produces Y_0, Y_1 \produces Y_0 X_0 \overline{Y}_0, X_0
  \produces a, Y_0 \produces b \} \\ 
\end{array} 
\right\rangle.$$
For example, $w(X_5) = baba\ovb \ova \ovb$.  Deduce that
$w(X_{k+2}) = w(Y_{k+1}) = ba \cdot \overline{w(Y_k)}$ for $k \geq
1$.
\end{example}


There is an important corollary of Theorem~\ref{Thm:Lohrey}:

\begin{theorem}[Lohrey~\cite{Lohrey04}]
\label{Thm:WordProbInFree}
The word problem for compressed words in the free group is solvable in
polynomial time. \qed
\end{theorem}


We notice two more consequences.

\begin{corollary}
\label{Cor:ComputeCyclicallyRed}
There is a polynomial-time algorithm that, given a compressed word
$\AA$, finds a compressed word $\XX$ with the property that $w_X$ is
the cyclic reduction of $w_A$.  The algorithm also gives the
compressed conjugating word. 
\end{corollary}

\begin{proof}
Following Theorems~\ref{Thm:Hagenah} and~\ref{Thm:Lohrey} assume that
$w_A$ is freely reduced.  Using the algorithm of
Lemma~\ref{Lem:ComputeInverse} produce the compressed word
$\overline{\AA}$.  Now apply the generalization of Plandowski's
algorithm (Theorem~\ref{Thm:PlandowskiII}) to find the largest $k$ so
that $w_A[:k] = w_{\overline{A}}[:k]$.  It follows that the
composition system $A' \produces A[:k]$ produces the conjugating word.
Also, the composition system $X \produces A[k:-k]$ produces the cyclic
reduction of $w_A$, as promised.
\end{proof}

The second consequence is more subtle:

\begin{theorem}
\label{Thm:ConjProbInFree}
The conjugacy problem for compressed words in the free group is
solvable in polynomial time.  The algorithm also computes the
compressed conjugating word.
\end{theorem}

We only sketch the proof, as the theorem is not used in the sequel.

\begin{proof}[Proof sketch of Theorem~\ref{Thm:ConjProbInFree}]
Suppose that $\AA$ and $\XX$ are the given compressed words.  Using
Corollary~\ref{Cor:ComputeCyclicallyRed} assume that $\AA$ and $\XX$
produce cyclically reduced words.  Using Lemma~\ref{Lem:ComputeLength}
check that $|w_A|$ and $|w_X|$ are equal.

Let $\WW$ be the compressed word with root production $W \produces X
\cdot X$.  That is $w_W = w_X \cdot w_X$.  Thus, to prove that $w_A$
and $w_X$ are conjugate it suffices to prove that $w_A$ appears as a
subword of $w_W$.  But this is exactly a special case of the {\em
fully compressed pattern matching problem} which can be solved in
polynomial time.  See, for example, the work of Karpinski, Rytter, and
Shinohara~\cite{KarpinskiEtAl95}, of Gasieniec \etal~\cite{GKPR96}, or
of Miyazaki, Shinohara, and Takeda~\cite{Miyazaki00}.
\end{proof}

\section{Free-by-cyclic groups}
\label{Sec:FreeByCyclic}

For group theory background the reader should consult Lyndon and
Schupp's book~\cite{LyndonSchupp77}.  Recall the definition of
$\Aut(F_m)$: the group of all automorphisms of the free group $F_m$.
Fix $\Phi \in \Aut(F_m)$.  The {\em free-by-cyclic} group $G_\Phi$ is
presented by:
$$\langle a_i, t \st ta_i\ovt = \Phi(a_i), i \in \{ 1, \ldots, m \}
\rangle.$$ 

The goal of this section is to prove:

\begin{theorem}
\label{Thm:FreeByCyclic}
The word problem for $G_\Phi$ is polynomial time.
\end{theorem}

This problem is already known to be in NP: Bridson and
Groves~\cite{BridsonGroves05} show that $G_\Phi$ has a quadratic
isoperimetric inequality.

\begin{proof}[Proof of Theorem~\ref{Thm:FreeByCyclic}]
Let $\calL_m = \{ a_i, \ova_i \}$ and $\calM = \calL_m \cup \{ t, \ovt
\}$.  Let $\calA = \{ A_{i,p} \st i \in \{ 1, \ldots m \}, p \in \NN
\}$.  Fix $\Phi \in \Aut(F_m)$ by assuming that the words $u_i(a_1,
\ldots, a_m) = \Phi(a_i)$ are given as input.

Define production rules as follows:
\begin{eqnarray*}
A_{i, 0} & \produces & a_i \\
A_{i, p} & \produces & u_i(A_{1, p-1}, \ldots, A_{m, p-1}), \quad p \geq 1
\end{eqnarray*}

Suppose now that $W$ is a word in $\calM^*$.  The length of $W$
determines the size of the given word problem.  Now rewrite $W$ in
stages: First freely reduce.  Next replace every $a_i$ and $\ova_i$
appearing by $A_{i,0}$ and by $\overline{A}_{i,0}$, respectively.  Now
move all occurrences of $t$ to the right, and of $\ovt$ to the left,
rewriting as follows:
\begin{eqnarray*}
t \cdot A_{i, p} & \leadsto & A_{i, p+1} \cdot t \\
t \cdot \overline{A}_{i, p} & \leadsto & 
                        \overline{A}_{i, p+1} \cdot t \\ 
A_{i, p} \cdot \ovt & \leadsto & \ovt \cdot A_{i, p+1} \\
\overline{A}_{i, p} \cdot \ovt & \leadsto & 
                     \ovt \cdot \overline{A}_{i, p+1} \\
t \cdot \ovt & \leadsto & \epsilon
\end{eqnarray*}


The result is a word in $\{A_{i,p}, \overline{A}_{i,p}\}^*$, possibly
with powers $t^k$ and $\ovt^l$ appearing at the end and beginning.
Let $W'(A_{i,p})$ be this word, omitting the leading and trailing
powers of $t$.  Construct a straight line program with a root
non-terminal $A \produces W'(A_{i,p})$.  Notice that $W$ is trivial in
$G_\Phi$ if and only if the word $w_A$ freely reduces in $F_m$ and the
powers satisfy $k = l$.  (This is a simple form of Britton's Lemma.
See page 181 of~\cite{LyndonSchupp77}.)

However the latter occurs if and only if the total exponent of $t$ in
$W$ is zero.  The former is exactly solved by applying Lohrey's
algorithm (Theorem~\ref{Thm:Lohrey}).
\end{proof}

\begin{remark}
\label{Rem:Acending}
The statement of Theorem~\ref{Thm:FreeByCyclic} may be generalized to
{\em ascending} HNN extensions: $\Phi$ is assumed to be an injection
instead of an automorphism.  The proof is identical.  The same
question for HNN extensions in general seems to be more delicate.
\end{remark}

\section{The automorphism group of a free group}
\label{Sec:Aut}

We now examine the automorphism group in greater detail.  Recall that
the automorphism group $\Aut(F_m)$ is finitely generated by the {\em
Nielsen generators} (see Chapter~1.4 of~\cite{LyndonSchupp77}):
\begin{enumerate}
\item
$\alpha_i \in \Aut(F_m)$ so that $\alpha_i|\calL_m$ interchanges $a_i$
  and $\ova_i$, fixing all other elements of $\calL_m$.
\item
$\beta_{ij} \in \Aut(F_m)$, with $i \neq j$, has $\beta_{ij}(a_i) =
  a_ia_j$, and $\beta_{ij}(a_k) = a_k$ for all $k \neq i$.
\end{enumerate}

\begin{remark}
\label{Rem:Unimportant}
Choosing a different generating set alters running times by at most a
multiplicative constant.  The choice above simplifies the proof below.
\end{remark}





\begin{theorem}
\label{Thm:Aut}
The word problem for $\Aut(F_m)$ is polynomial time.
\end{theorem}

This solves problem (C1) on the list maintained by
Baumslag, Myasnikov, and Shpilrain~\cite{Baumslag02}. 

\begin{proof}[Proof of Theorem~\ref{Thm:Aut}]
Suppose that $\Phi = \varphi_1 \ldots \varphi_n$ is a word in the
Nielsen generators of $\Aut(F_m)$.  We must check that $\Phi(a_i)$
freely reduces to $a_i$, for all $i$.

To do this, define a straight line program: Let $\{ A_{i, p} \}$ be
the set of non-terminals with $i \in \{ 1, \ldots m \}$ and $p \in \{
0, 1, \ldots n \}$.  Create the following production rules:
\begin{eqnarray*}
A_{i, 0} & \produces & a_i \\
A_{i, p} & \produces & \varphi_p(A_{i, p-1}), \quad p \geq 1
\end{eqnarray*}
If $\varphi_p = \alpha_i$ is of the first kind then $\varphi_p(A_{i,
p})$ equals $\overline{A}_{i, p-1}$ while $\varphi_p(A_{j, p})$ equals
$A_{j, p-1}$, for $j \neq i$.  If $\varphi_p = \beta_{ij}$ is of the
second kind then $\varphi_p(A_{i, p})$ equals $A_{i,p-1} \cdot A_{j,
p-1}$, and so on.

Now apply the algorithm of Theorem~\ref{Thm:Lohrey} to rewrite this
straight line program so that all outputs are freely reduced.  If the
resulting composition system has $|w(A_{i, n})| \geq 2$ for any $i$
then the automorphism $\Phi$ is nontrivial.  If $|w(A_{i, n})| = 1$
for all $i$ then, using the algorithm of
Lemma~\ref{Lem:ComputeSubword}, check if $w(A_{i,n}) = a_i$.  If this
is the case for all $i$ then $\Phi$ is the identity element of
$\Aut(F_m)$.
\end{proof}

\begin{remark}
\label{Rem:FreeByCyclicRevisit}
In our analysis of the word problem for free-by-cyclic groups we could
accept as input both the word $W$ in $\calM^*$ and automorphism $\Phi$
given as a word in the Nielsen generators.  Now we need not precompute
the words $u_i$.  Instead we find, as in the proof of
Theorem~\ref{Thm:Aut}, straight line programs producing these words.
The running time of Theorem~\ref{Thm:FreeByCyclic} then becomes
polynomial in the two inputs $W$ and $\Phi$.
\end{remark}

Note that solving the word problem for a group also solves the word
problem for subgroups.  Of course, there are many beautiful subgroups
of $\Aut(F_m)$.  As just a single example consider the {\em braid
group}, $B_m$: Let $\DD_m$ be a disk with $m$ points removed from the
interior.  Then $B_m$ is the group of homeomorphisms of $\DD_m$ which
fix the boundary pointwise, modulo boundary and puncture fixing
isotopies.  Choosing a basepoint on the boundary makes $B_n$ act on
$\pi_1(\DD_m) \isom F_m$ and so embeds $B_n$ into $\Aut(F_m)$.  A
simple corollary to Theorem~\ref{Thm:Aut} is the well-known:

\begin{corollary}[\cite{Bigelow01}, \cite{Krammer02},
    \cite{EpsteinEtAl92}, \cite{BirmanEtAl98},
    \cite{Hamidi-Tehrani00}]
\label{Cor:Braid}
The word problem for $B_m$ is polynomial time. \qed
\end{corollary}


\section{Membership problems}
\label{Sec:Membership}

We now turn our attention to {\em membership problems}, also called
{\em generalized word problems}.  Suppose that $H$ is a subgroup of a
finitely presented group $G$.  We seek an algorithm that, given a word
$W$ written in the generators of $G$, decides if $W$ represents an
element of $H$.  Note that if $H$ is normal in $G$ then such an
algorithm also solves the word problem in the quotient $G/H$.


Recall that $\Inn(F_m)$ is the normal subgroup of $\Aut(F_m)$ which
contains $\Phi_U$ for all words $U \in \calL_m^*$: $\Phi_U(V) =
UV\overline{U}$.  The quotient is the {\em outer automorphism} group,
$\Out(F_m)$.

\begin{theorem}
\label{Thm:Out}
The word problem for $\Out(F_m)$ is polynomial time. 
\end{theorem}

\begin{proof}
To see this, note that the membership problem for $\Inn(F_m)$, inside
of $\Aut(F_n)$, is solved by the construction given in the proof of
Theorem~\ref{Thm:Aut} and by Corollary~\ref{Cor:ComputeCyclicallyRed}.
\end{proof}

Another interesting membership problem is that of the braid group (or
more generally, {\em mapping class groups} of punctured surfaces)
inside of $\Aut(F_m)$.  In order to avoid multiplying examples we only
discuss the problem for the braid group: it suffices to check that the
boundary word is fixed and all punctures are preserved, up to
conjugacy.  (See~\cite[Theorem~N6]{MagnusEtAl66}, for example.)
Corollary~\ref{Cor:ComputeCyclicallyRed} can check the latter and the
former is dealt with by Lohrey's Algorithm (Theorem~\ref{Thm:Lohrey}).

Here is another kind of membership problem, first proposed by Nathan
Broaddus.  Let $\MCG(S_g)$ denote the mapping class group of the
closed connected orientable genus $g$ surface, $S_g$.  Consider a
handlebody $V_g$ so that $\bdy V_g = S_g$.  Then $\MCG(V)$ naturally
includes in $\MCG(S)$.

\begin{remark}
\label{Rem:Meridians}
Recall the fundamental fact that $\Phi \in \MCG(S)$ lies in the
subgroup $\MCG(V)$ if and only if $\Phi$ preserves the set of {\em
meridians}: the set of curves in $S$ which bound disks in $V$.  In
fact a ``weaker'' condition is equivalent: let $\DD$ be a collection
of $g$ disjoint disks in $V$ so that $V \setminus \DD$ is a
three-ball. Then $\Phi \in \MCG(V)$ if and only if the curves
$\Phi(\bdy \DD)$ bound disks in $V$.
\end{remark}

\begin{figure}[ht]
\psfrag{}{}
\psfrag{}{}
$$\begin{array}{c}
\epsfig{file=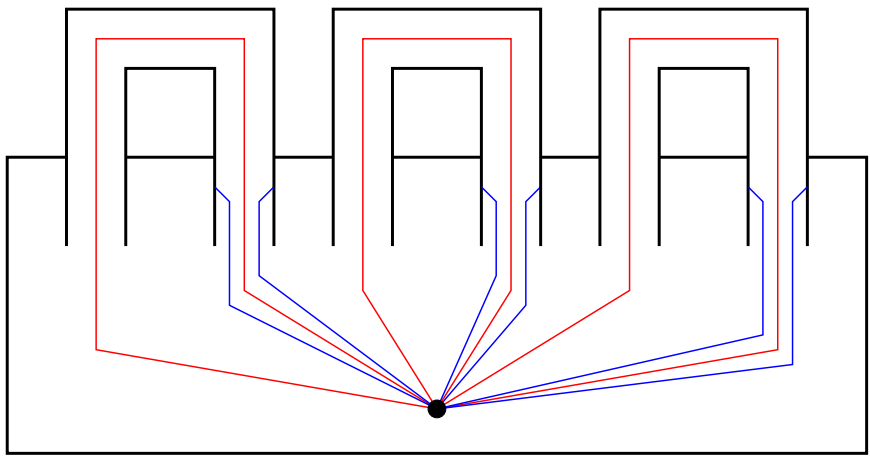, height = 6 cm}
\end{array}$$
\caption{Generators for the fundamental group of $S$.}
\label{Fig:Surface}
\end{figure}

Fix a point $x \in S$ and let $\{ a_1, \ldots, a_g, b_1, \ldots, b_g
\}$ be the standard set of generators of $\pi_1(S, x)$.  We arrange
matters so that all of the $b_i$ are meridians. See
Figure~\ref{Fig:Surface}.

\begin{figure}[ht]
\psfrag{}{}
\psfrag{}{}
$$\begin{array}{c}
\epsfig{file=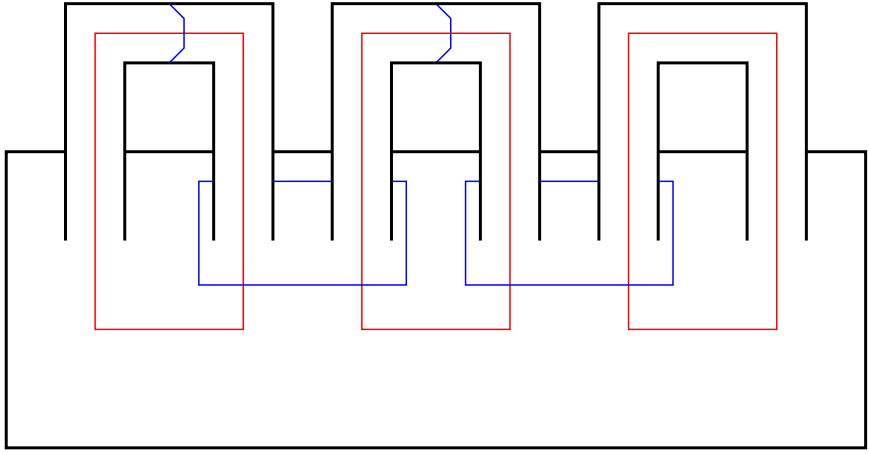, height = 6 cm}
\end{array}$$
\caption{Generators for the mapping class group.}
\label{Fig:Twist}
\end{figure}

Choose also a standard set $\{ \tau \}$ of Dehn twist generators for
$\MCG(S)$. In fact we will over-specify these twists: for each Dehn
twist in the generating set choose a twisting curve $\{ \alpha \}$
which avoids a small neighborhood of the basepoint $x \in S$.  See
Figure~\ref{Fig:Twist}.

Pick finally a point $y$ in this small neighborhood of $x$, avoiding
all the loops $a_i$ and $b_i$.  It follows that $\pi_1(S \setminus \{
y \}, x) \isom F_{2g}$ is freely generated by the $a_i$ and $b_i$
loops.  Also, the Dehn twists give free group automorphisms.
(See~\cite[Theorem~N10]{MagnusEtAl66}.)
Broaddus tells us the useful:

\begin{lemma}
\label{Lem:Broaddus}
Fix $W \in \{ a_i, b_i, \ova_i, \ovb_i \}^*$ so that $W$ is homotopic
in $S$ to a simple closed curve $\omega$.  The following two
conditions are equivalent:
\begin{itemize}
\item
$\omega$ is a meridian.  
\item
After deleting all $b_i$'s and $\ovb_i$'s from $W$ the resulting word
in $\{ a_i, \ova_i \}^*$ freely reduces to the empty word.
\end{itemize}
\end{lemma}

\begin{proof}
Consider the map $\pi_1(S,x) \to \pi_1(V,x) \isom F_g$ induced by
inclusion.  As in \reffig{Surface} all of the $b_i$ lie in the kernel
while the images of the $a_i$ freely generate $\pi_1(V)$.
Accordingly, identify $\pi_1(V)$ and the free group $\langle a_i \st
\rangle$.

Suppose now that the first condition holds: $\omega$ is a meridian.
Then $W$ lies in the kernel.  Since the $b_i$ normally generate the
kernel the second condition follows.

Suppose instead that the second condition holds.  Since $W$ is in the
kernel, deduce that $\omega$ bounds a singular disk in $V$.  By the
celebrated Loop Theorem~\cite{Hempel76} $\omega$ is a meridian.
\end{proof}

We are now prepared to prove:

\begin{theorem}
\label{Thm:Handlebody}
The membership problem for $\MCG(V)$ in $\MCG(S)$ is polynomial time.
\end{theorem}

\begin{proof}
As is well-known~\cite{Humphries79} the $2g + 1$ Dehn twists shown in
Figure~\ref{Fig:Twist} generate the mapping class group, $\MCG(S)$.
So fix a word $\Phi = \tau_1\ldots\tau_n$, written in terms of these
twists and their inverses.  By Remark~\ref{Rem:Meridians} it is enough
to check that $\Phi(b_i)$, thought of as a word in the free group
generated by the $a_i$ and the $b_i$, satisfies the second condition
of Lemma~\ref{Lem:Broaddus}.  This could be done directly, but
$\Phi(b_i)$ might have length exponential in $n$.

Instead, for each $i$ encode $\Phi(b_i)$ as a compressed word, say
$\AA(i)$.  It is a triviality to remove all $b_i$ and $\ovb_i$
appearing in $w(A(i))$: for every non-terminal $B$ with $B \produces
b_i$ replace the production by $B \produces \epsilon$.  Call the new
compressed word $\AA'(i)$.

Now run Lohrey's algorithm (Theorem~\ref{Thm:Lohrey}) on $\AA'(i)$.  The
mapping class $\Phi$ lies in $\MCG(V)$ if and only if $w(A'(i))$ freely
reduces to the empty string, for all $i$.
\end{proof}


\begin{remark}
\label{Rem:Heegaard}
Suppose that $W$ is another handlebody of genus $g$ with $\bdy W = S$
so that $V \cup_S W$ is the three-sphere. Here the $a_i$ loops bound
disks in $W$.  Recall that the {\em Heegaard group} $\calH = \MCG(V)
\cap \MCG(W)$ is the subgroup of $\MCG(S)$ consisting of all mapping
classes which extends over both $V$ and $W$.  Applying
Theorem~\ref{Thm:Handlebody} twice gives a polynomial-time algorithm
for the membership problem of $\calH$ in $\MCG(S)$.
\end{remark}



\section{Hagenah's algorithm}
\label{Sec:Hagenah}

We now discuss the computer science underpinnings of the discussion
above.  To begin, in Hagenah's thesis~\cite[Chapter~8]{Hagenah00} we
find:

\begin{theorem}
\label{Thm:Hagenah}
There is a polynomial-time algorithm that, given a composition system
$\AA$, finds a straight line program $\XX$ with the property that $w_X
= w_A$.
\end{theorem}

The exposition of this result in~\cite{Hagenah00} is wonderfully
clear.  I present a proof only to make this paper self-contained.

\begin{proof}[Proof of Theorem~\ref{Thm:Hagenah}]
Fixing notation, suppose that $\AA = \langle \calL, \calA, A, \calP
\rangle$.  Note that all productions in $\calP$ are either of the form
$B \produces C[i:j] \cdot D[k:l]$ or of the form $B \produces a[i:j]$.

Build the straight line program $\XX = \langle \calL, \calX, X, \calQ
\rangle$ from the bottom up.  The set $\calX$ will contain {\em plain}
non-terminals, one for each non-terminal of $\calA$, and {\em
decorated} non-terminals, each associated to some plain non-terminal.
Proceed as follows: for every non-terminal $B \in \calA$ of height one
we add a plain non-terminal $Y$ to $\calX$.  Suppose $B$ produces
$a[i:j]$. Then 
\begin{itemize}
\item
if $j = i+1$ add $Y \produces a$ to $\calQ$ and
\item
if $j = i$ add $Y \produces \epsilon$ to $\calQ$.
\end{itemize}

Let $A$ be the root non-terminal of $\calA$.  Assume via induction
that for every other non-terminal $B \in \calA$ a plain non-terminal
$Y$ has been added to $\calX$, so that $w_Y = w_B$.  We now describe
the decorated non-terminals that may also, by induction, appear in
$\calX$.  Fix any plain non-terminal $Y$ in $\calX$.  Then $Y^{[i:j]}$
is a {\em decorated} non-terminal.  There are various cases:
\begin{itemize}
\item
If $0 < i < j < |w_Y|$ then $Y^{[i:j]}$ is a {\em subword} non-terminal.
\item
If $0 < j < |w_Y|$ then $Y^{[:j]}$ is a {\em prefix} non-terminal.
\item
If $0 < i < |w_Y|$ then $Y^{[i:]}$ is a {\em suffix} non-terminal.
\item
$Y^{[:]} = Y$ is the plain non-terminal.
\item
$Y^{[i:i]} = \epsilon$ is the empty word. 
\end{itemize}
Repeated decoration behaves as expected: $(Y^{[i:j]})^{[l:k]} =
Y^{[i+k:i+l]}$.  The production rules for decorated non-terminals are
given below.

Suppose now that the root has production $A \produces B[i:j] \cdot
C[k:l]$.  Suppose that $Y$ and $Z$ are plain non-terminals in $\calX$
corresponding to the non-terminals $B$ and $C$.  Add a plain
non-terminal $X$ to $\calX$ corresponding to $A$.  Add the
non-terminals $Y^{[i:j]}$ and $Z^{[k:l]}$ to $\calX$.  Add the
production rule $X \produces Y^{[i:j]} \cdot Z^{[k:l]}$ to $\calQ$.  

Production rules are needed for every new decorated non-terminal,
$Y^{[i:j]}$, created by the addition of the plain $X$.  Suppose that
the plain non-terminal $Y$ produces $U \cdot V$.  (Here it is possible
that $U$ and $V$ are themselves decorated non-terminals.)  There are
several cases and subcases: Suppose first that $Y^{[i:j]}$ is a
subword non-terminal.
$$\begin{array}{lll}
\mbox{Subcase}      &  \mbox{Add to $\calX$}  &  \mbox{Add to $\calQ$} \\
                                                          \hline
|w_U| \leq i        & V^{[i-|w_U|:j-|w_U|]}          & Y^{[i:j]} \produces V^{[i-|w_U|:j-|w_U|]} \\
i < |w_U| < j \quad & U^{[i:]}, V^{[:j-|w_U|]} \quad & Y^{[i:j]} \produces U^{[i:]} \cdot V^{[:j-|w_U|]} \\
j \leq |w_U|        & U^{[i:j]}                      & Y^{[i:j]} \produces U^{[i:j]}\\
\end{array}$$

\noindent
Suppose now that $Y^{[:j]}$ is a prefix non-terminal.
$$\begin{array}{lll}
\mbox{Subcase}       &  \mbox{Add to $\calX$}  \quad &  \mbox{Add to $\calQ$} \\
                                                          \hline
|w_U| < j \quad & V^{[:j-|w_U|]} & Y^{[:j]} \produces U \cdot V^{[:j-|w_U|]} \\
j \leq |w_U|    & U^{[:j]}       & Y^{[:j]} \produces U^{[:j]}
\end{array}$$

\noindent
Suppose now that $Y^{[i:]}$ is a suffix non-terminal.
$$\begin{array}{lll}
\mbox{Subcase}       &  \mbox{Add to $\calX$}  \quad &  \mbox{Add to $\calQ$} \\
                                                          \hline
|w_U| \leq i \quad & V^{[i-|w_U|:]} & Y^{[i:]} \produces V^{[i-|w_U|:]} \\
i < |w_U|          & U^{[i:]}       & Y^{[i:]} \produces U^{[i:]} \cdot V 
\end{array}$$

Notice that creating the plain non-terminal $X$ causes at most two
decorated non-terminals to be created, both of lesser height.  Every
subword non-terminal in turn creates at most one subword non-terminal
or at most one prefix and at most one suffix non-terminal.  Again,
these have lesser height.  Finally, any prefix (suffix) non-terminal
causes at most one prefix (suffix) non-terminal to be created.  As
usual, the height decreases.

Suppose that $n = ||A||$.  It follows that the creation of the plain
non-terminal $X$ adds at most $1 + 2(2n) = 1 + 4n$ new decorated
non-terminals to $\calX$.  Thus the total number of non-terminals in
$\calX$, at the end of the construction, is 
$$n + \sum_{B \in \calA} (1 + 4||B||) \leq n + \sum_{i = 1}^n (1 + 4i)
   = 2n^2 + 4n.$$
This completes both the description of the algorithm and its
proof of correctness.
\end{proof}

\section{Plandowski's Algorithm}
\label{Sec:Plandowski}

The final piece of the puzzle is:

\begin{theorem}[Plandowski~\cite{Plandowski94}]
\label{Thm:PlandowskiI}
There is a polynomial-time algorithm that, given straight line
programs $\AA$ and $\XX$ in normal form, decides whether or not $w_A =
w_X$.
\end{theorem}

A proof, essentially following~\cite{Plandowski94}, is provided for
the convenience of the reader.

\begin{proof}[Proof of Theorem~\ref{Thm:PlandowskiI}]
Let $\AA = \langle \calL, \calA, A, \calP \rangle$ and $\XX = \langle
\calL, \calX, X, \calQ \rangle$.  Note that we assume, as we may, that
$\AA$ and $\XX$ have the same terminal alphabet.  Making a copy of
$\XX$ if necessary, assume that $\calA \cap \calX = \emptyset$.
Finally assume that $|w_A| = |w_X|$ and $|\calA| = m \geq n =
|\calX|$.

We begin with the following definition: a triple $(B, Y, i)$ is an
{\em assertion} if:
\begin{itemize}
\item
$B \in \calA \cup \calL$ and $Y \in \calX \cup \calL$.
\item
$0 \leq i < |w_B|$.
\end{itemize}
If $0 \leq i$ and $|w_B| \leq i + |w_Y|$ then $(B, Y, i)$ is a {\em
overlap} assertion.  If $0 < i$ and $i + |w_Y| < |w_B|$, and then $(B,
Y, i)$ is a {\em subword} assertion.  Assertions of the form $(Y, B,
i)$, are defined similarly.  We do not allow a pair of non-terminals
from the same program to appear in a single assertion.

An overlap assertion $(B, Y, i)$ is {\em satisfied} if and only if
$w_B[i:] = w_Y[:|w_B|-i]$.  Likewise, a subword assertion is satisfied
if and only if $w_B[i:i+|w_Y|] = w_Y$.  As a bit of terminology a set
of assertions, $\Gamma$, is satisfied if and only if every assertion
$\gamma \in \Gamma$ is.  In point of fact, the algorithm checks
satisfiability of $(B, Y, i)$ when and only when both $B$ and $Y$ are
terminal characters and $i = 0$.

In general a set of assertions $\Gamma_k$ is transformed into another
such set, $\Gamma_{k+1}$.  Beginning with $\Gamma_0 = \{ (A, X, 0) \}$
the following properties will be maintained:
\begin{enumerate}
\item[(a)]
$\Gamma_{k+1}$ is satisfied if and only if $\Gamma_{k}$ is satisfied.
\item[(b)]
At most $m+n-k$ elements of $\calA \cup \calX$ are mentioned in
$\Gamma_k$.
\item[(c)] 
For all $k$, $|\Gamma_k|$ is bounded by $(k+1)4mn(m+n)$.
\end{enumerate}
There are two ways to produce a new assertions from old, {\em
splitting} and {\em compacting}.  

\subsection*{Splitting}

Fix $\Gamma$, a set of assertions. Suppose that, of all non-terminals
from $\calA$ and $\calX$ appearing in $\Gamma$, the non-terminal $B
\in \calA$ has maximal length.  Fix $\gamma \in \Gamma$.  We must
define $\Split(\gamma, B)$ and then $\Split(\Gamma, B)$.

There are several cases to consider.  If $B$ does not appear in
$\gamma$ then $\Split(\gamma, B) = \{ \gamma \}$.  Now suppose that
$\gamma = (B, Y, i)$ or $\gamma = (Y, B, i)$.  Note that $\gamma$ is
either an overlap or subword assertion and that we have assumed $|w_B|
\geq |w_Y|$.  Suppose that $B \produces C \cdot D$.  Now consider
subcases.  If $\gamma = (B, Y, i)$ is an overlap assertion then:
$$\begin{array}{lll}
\mbox{Subcase}  &  \Split(\gamma, B)  &  \mbox{Type} \\
                                                          \hline
i < |w_C|          & (C, Y, i)             & \mbox{overlap} \\
                   & (Y, D, |w_c|-i)       & \mbox{either}  \\ 
|w_C| \leq i \quad & (D, Y, i-|w_c|) \quad & \mbox{overlap} 
\end{array}$$

\noindent
The table should be read as follows: When $i < |w_C|$ then
$\Split(\gamma, B)$ contains two assertions: either a pair of overlaps
or one of each type.  When $|w_C| \leq i$ the set $\Split(\gamma, B)$
contains a single overlap assertion.  Suppose now that $\gamma = (B,
Y, i)$ is a subword assertion:
$$\begin{array}{lll}
\mbox{Subcase}  &  \Split(\gamma, B)  &  \mbox{Type} \\
                                                          \hline
i+|w_Y| < |w_C|            & (C, Y, i)             & \mbox{subword} \\
i+|w_Y| = |w_C|            & (C, Y, i)             & \mbox{overlap} \\
i < |w_C| < i+|w_Y|  \quad & (C, Y, i)             & \mbox{overlap} \\
                           & (Y, D, |w_c|-i) \quad & \mbox{overlap} \\  
i = |w_C|                  & (Y, D, 0)             & \mbox{overlap} \\
|w_C| < i                  & (D, Y, i-|w_C|)       & \mbox{subword} 
\end{array}$$

\noindent
Suppose now that $\gamma = (Y, B, i)$ is an overlap assertion.  As
usual assume that $|w_B| \geq |w_Y|$:
$$\begin{array}{lll}
\mbox{Subcase}  &  \Split(\gamma, B)  &  \mbox{Type} \\
                                                          \hline
|w_Y| \leq i+|w_c| \quad & (Y, C, i)             & \mbox{overlap} \\
i+|w_C| < |w_Y|          & (Y, C, i)             & \mbox{subword} \\
                         & (Y, D, i+|w_c|) \quad & \mbox{overlap} \\
\end{array}$$

\noindent
Finally $(Y, B, i)$ cannot be a subword assertion because $|w_B| \geq
|w_Y|$.  This finishes the definition of $\Split(\gamma, B)$.  Define
$$\Split(\Gamma, B) = \cup_{\gamma \in \Gamma} \Split(\gamma, B).$$
Immediate from the definitions is:

\begin{claim}
\label{Clm:SplitSat}
A set of assertions $\Gamma$ is satisfied if and only if
$\Split(\Gamma, B)$ is satisfied. \qed
\end{claim}

Define now $o(\Gamma)$ to be the number of overlap assertions in
$\Gamma$.  Similarly we take $s(\Gamma)$ to be the number of subword
assertions of $\Gamma$.  So $|\Gamma| = o(\Gamma) + s(\Gamma)$.  From
the tables above deduce:

\begin{claim}
\label{Clm:SplitSize}
Suppose that $\Gamma$ is a set of assertions. Then
\begin{eqnarray*}
o(\Split(\Gamma, P)) & \leq & o(\Gamma) + 2s(\Gamma) \\
s(\Split(\Gamma, P)) & \leq & o(\Gamma) + s(\Gamma) 
\end{eqnarray*} 
when $P \in \calA \cup \calX$ is a non-terminal of maximal length
in $\Gamma$.  \qed
\end{claim}

\subsection*{Compact}

Now for the definition of $\Compact(\Gamma)$.  Note that, if $u$ is a
word, then $p \in \NN$ is a {\em period} of $u$ if
\begin{itemize}
\item
$1 \leq p \leq |u| - 1$ and 
\item
$u[i] = u[i + p]$ for all $0 \leq i \leq |u| - 1 - p$.
\end{itemize}
An immediate consequence of the definition is:

\begin{claim}
\label{Clm:EasyPeriod}
Suppose that $\gamma = (B, Y, i)$ and $\gamma' = (B, Y, j)$ are
overlap assertions with $i < j$.  Then $\gamma$ and $\gamma'$ are
satisfied if and only if $\gamma$ is satisfied and $j - i$ is a period
of the word $w_B[i:]$. \qed
\end{claim}

We now give a restricted version of the famous {\em Periodicity
Lemma}~\cite{Lothaire83}:

\begin{lemma}
\label{Lem:Period}
If $p$ and $q$ are periods of $u$, where $p + q \leq |u|$, then
$\gcd(p, q)$ is also a period of $u$. \qed
\end{lemma}


The following claim is the engine in the proof of correctness of
Plandowski's Algorithm:

\begin{claim}
\label{Clm:SubtlePeriod}
Suppose that $\gamma = (B, Y, i)$, $\gamma' = (B, Y, j)$, and
$\gamma'' = (B, Y, k)$ are overlap assertions with $i < j < k$ and $j
- i + k - i \leq |w_B| - i$.  Then $\gamma$, $\gamma'$, and $\gamma''$
are satisfied if and only if $\gamma$ is satisfied and $\gcd(j - i, k
- i)$ is a period of $w_B[i:]$.
\end{claim}

\begin{proof}
This follows from two applications of Claim~\ref{Clm:EasyPeriod} and
from the Periodicity Lemma~\ref{Lem:Period}.
\end{proof}

Equivalently $\gamma$, $\gamma'$ and $\gamma''$ are satisfied if and
only if $\gamma$ and $\delta = (B, Y, i + \gcd(j - i, k - i))$ are
satisfied.  This leads directly to the definition of $\SimpleCompact$:
given $\{ \gamma, \gamma', \gamma'' \}$ as in the hypothesis of
Claim~\ref{Clm:SubtlePeriod} define $\SimpleCompact(\{ \gamma,
\gamma', \gamma'' \}) = \{ \gamma, \delta \}$, with $\delta$ as above.

Now, for any set of assertions $\Gamma$ define $\Compact(\Gamma)$ to
be the result of applying $\SimpleCompact$ to all triples of overlap
assertions which follow the requirements of
Claim~\ref{Clm:SubtlePeriod}.  Every successful application of
$\SimpleCompact$ removes an assertion from $\Gamma$.  Thus a single
$\Compact$ operation involves calling $\SimpleCompact$ at most
$O(|\Gamma|^4)$ times. (This can be greatly improved upon, if so
desired.)  We also record the fact:

\begin{claim}
\label{Clm:CompactSat}
$\Gamma$ is satisfied if and only if $\Compact(\Gamma)$ is
satisfied. \qed
\end{claim}

Now to define $\Gamma_{k+1}$ in terms of $\Gamma_k$.  Suppose that $P
\in \calA \cup \calX$ is a non-terminal appearing in $\Gamma_k$
maximizing the length of $|w_P|$.  Then take
$$\Gamma_{k+1} = \Compact(\Split(\Gamma_k, P)).$$ Note that property
(a) above is guaranteed by Claims~\ref{Clm:SplitSat}
and~\ref{Clm:CompactSat} while property (b) is provided by the fact
that every non-terminal is split for at most one value of $k$.  We
must now bound the size of $\Gamma_k$.

Fix attention on any pair of non-terminals $B \in \calA$ and $Y \in
\calX$.  Let $\{(B, Y, i_j)\}_{j = 1}^N$ be the overlap assertions of
$\Gamma_k$ which mention $B$ and $Y$ in that order and indexed so that
$i_j < i_{j+1}$.  As $\Gamma_k$ is compact it follows from
Claim~\ref{Clm:SubtlePeriod} that
$$|w_B| - i_j < i_{j+2} - i_j + i_{j+1} - i_j.$$ Since $i_{j+1} <
i_{j+2}$ it follows that
$$\frac{1}{2}(|w_B| - i_j) < i_{j+2} - i_j.$$ Deduce that $N \leq 2
\log_2(|w_B|) + 1 \leq 2||B|| + 1$, with the last inequality following
from Lemma~\ref{Lem:LengthBound}.

\begin{claim}
\label{Clm:CompactSize}
With $\Gamma_{k+1}$ as given:
\begin{eqnarray*}
o(\Gamma_{k+1}) & \leq & mn(2m+1) + nm(2n+1) \leq 4mn(m+n) \\
s(\Gamma_{k+1}) & \leq & 4mn(m+n) + s(\Gamma_k) \leq k(4mn(m+n)).
\end{eqnarray*} 
This verifies property (c) above. \qed
\end{claim}

This completes both the description of the algorithm and its
proof of correctness.
\end{proof}

\appendix
\section{On surfaces}
\label{Sec:Surfaces}

The discussion above gives a satisfactory picture of the behavior of
compressed words in the free group.  One immediately asks for a
similar treatment of hyperbolic groups in general.  However the
situation there appears to require a new idea.

Instead we briefly describe {\em well-tempered paths}: a beautiful
geodesic language for closed surface groups.  Our discussion is meant
to be more inspiring than exhaustive: many details are omitted.  For
simplicity, we restrict ourselves to $S = S_2$ the closed, orientable,
connected genus two surface.

\begin{figure}[ht]
\psfrag{a}{$a$}
\psfrag{b}{$b$}
\psfrag{c}{$c$}
\psfrag{d}{$d$}
\psfrag{e}{$e$}
$$\begin{array}{c}
\epsfig{file=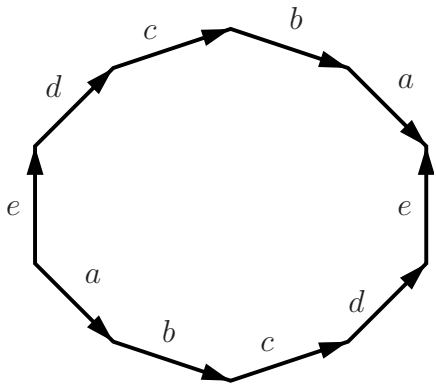, height = 5 cm}
\end{array}$$
\caption{A labelled regular hyperbolic decagon.}
\label{Fig:Decagon}
\end{figure}

Let $D$ be the regular decagon in the hyperbolic plane with angles
$2\pi/5$.  Label the boundary of $D$ with the word $abcdeabcde$, read
counter-clockwise.  The first five edges are oriented
counter-clockwise while the last five are oriented clockwise. Let
$\calL = \{ a, b, c, d, e, \ova, \ovb, \ovc, \ovd, \ove \}$.  See
Figure~\ref{Fig:Decagon}.

The decagon and its labelling extends to a tiling $\calD$ of the
hyperbolic plane.  Notice that a path in the one-skeleton determines a
word in $\calL^*$.  Let $\calM \subset \calL^*$ be the subset that can
be realized in this way.  Conversely, a word of $\calM$ determines a
unique path in the tiling, up to the action of $\pi_1(S)$.  When it
cannot cause confusion we treat paths and words interchangeably.

\begin{figure}[ht]
\psfrag{a}{$a$}
\psfrag{b}{$b$}
\psfrag{c}{$c$}
\psfrag{d}{$d$}
\psfrag{e}{$e$}
$$\begin{array}{c}
\epsfig{file=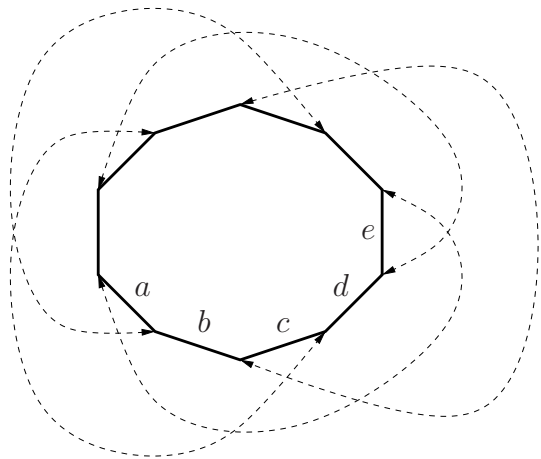, height = 6 cm}
\end{array}$$
\caption{All of the good turns are shown.}
\label{Fig:Turns}
\end{figure}

The following words of length five are called {\em bad turns}:
$$edcba, ~\ovb\ova edc, ~\ovd\ovc\ovb\ova e, ~a\ove\ovd\ovc\ovb,
~cba\ove\ovd,$$
$$\ova\ovb\ovc\ovd\ove, ~\ovc\ovd\ove ab, ~\ove abcd, ~bcde\ova,
~de\ova\ovb\ovc.$$ A path is {\em well-tempered} if the corresponding
word is freely reduced and contains no bad turn.  The intent, and
hence the name, is that these paths want to ``turn right'' as often as
they ``turn left.''  This is possible because $10/2$ is odd.  See
Figure~\ref{Fig:Turns} for a picture of the good turns.

Notice that well-tempered paths are well-behaved:

\begin{theorem}
\label{Thm:WellTempered}
For any ordered pair of vertices in the tiling $\calD$ there is a
unique well-tempered path connecting one to the other.  Well-tempered
paths are geodesic.  Subpaths are again well-tempered as are inverses.
Finally, well-tempered paths are {\em locally detectable}: to verify
the property it suffices to check all subpaths of length five.
\end{theorem}


\begin{remark}
\label{Rem:ShortLex}
Well-tempered words are similar to {\em short-lex} words in hyperbolic
groups (see~\cite{EpsteinEtAl92}).  For example, both form regular
languages which satisfy uniqueness and which are closed under taking
subwords.  Short-lex is generally not closed under taking inverses and
is essentially never locally detectable.
\end{remark}


Here is a sketch of the proof of Theorem~\ref{Thm:WellTempered}.
Suppose that $\alpha$ and $\beta$ are paths and the final vertex of
$\alpha$ is the initial vertex of $\beta$.  Define $\alpha
\cdot \beta$ to be their concatenation.  If $\alpha$ and $\beta$ are
well-tempered and $\alpha \cdot \beta$ is not then the free reduction
(or bad turns) must overlap the point of concatenation.

We may straighten $\alpha \cdot \beta$ (rel endpoints) until it
becomes well-tempered.  There are four stages:
\begin{enumerate}
\item
Free reduction.
\item
Sweeping across two sides of at most one corridor.
\item
Sweeping across at most three pieces.
\item
Sweeping across one side of at most two corridors.
\end{enumerate}

\begin{figure}[ht]
$$\begin{array}{c}
\psfrag{e}{$e$}
\epsfig{file=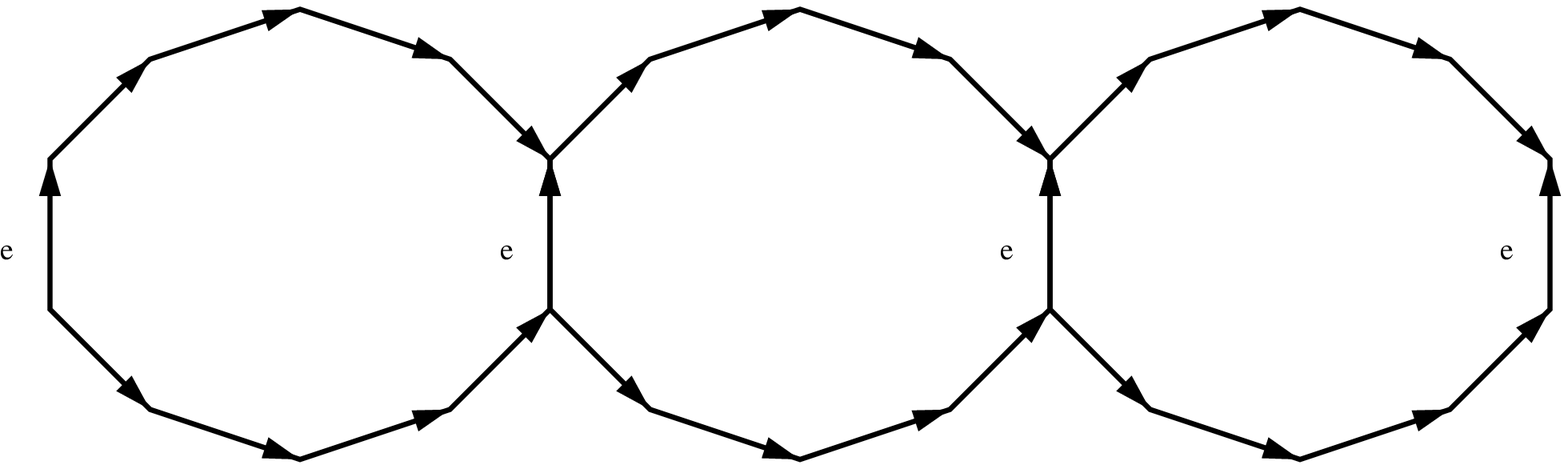, height = 3 cm}
\end{array}$$
\caption{An $e$-corridor of length three.}
\label{Fig:Corridor}
\end{figure}

Before elaborating on these we briefly give definitions: A {\em piece}
is a path $\delta$ of length two to ten where all edges of $\delta$
are on the boundary of a single decagon of the tiling.  Now, fix a
decagon $D = D_0$ and suppose that $D_i$ is the image of $D_0$ under
the $i^\th$ power of a fixed side pairing transformation of $D_0$.
The union $C = \cup_{i=0}^{k-1} D_i$ is called a {\em corridor}.  The
two edges in $\bdy C$ corresponding to the transformation are the {\em
ends} of $C$.  The other two components of $\bdy C$ are the {\em
sides} of $C$.  Note that sides of corridors always have period four.
See Figure~\ref{Fig:Corridor}.

Suppose again that $\alpha$ and $\beta$ are well-tempered and
$\gamma_0 = \alpha \cdot \beta$ is the concatenation.  We can now
flesh out the stages required to make $\gamma_0$ well-tempered,
assuming it is not already. First freely reduce, if possible, to
produce $\gamma_1$.  If $\gamma_1$ contains two sides and one end of a
corridor $C$ then sweep $\gamma_1$ across $C$ to obtain $\gamma_2$.
This deals with all pieces of length nine and ten.

Next sweep $\gamma_2$ across at most three pieces of lengths between
five and eight to form $\gamma_3$.  The pieces of length five are
necessarily bad turns.  The proof that there are at most three such is
a lengthy but straight-forward combinatorial argument.

\begin{figure}[ht]
\psfrag{a}{$\alpha$}
\psfrag{b}{$\beta$}
\psfrag{c}{$\gamma$}
$$\begin{array}{c}
\epsfig{file=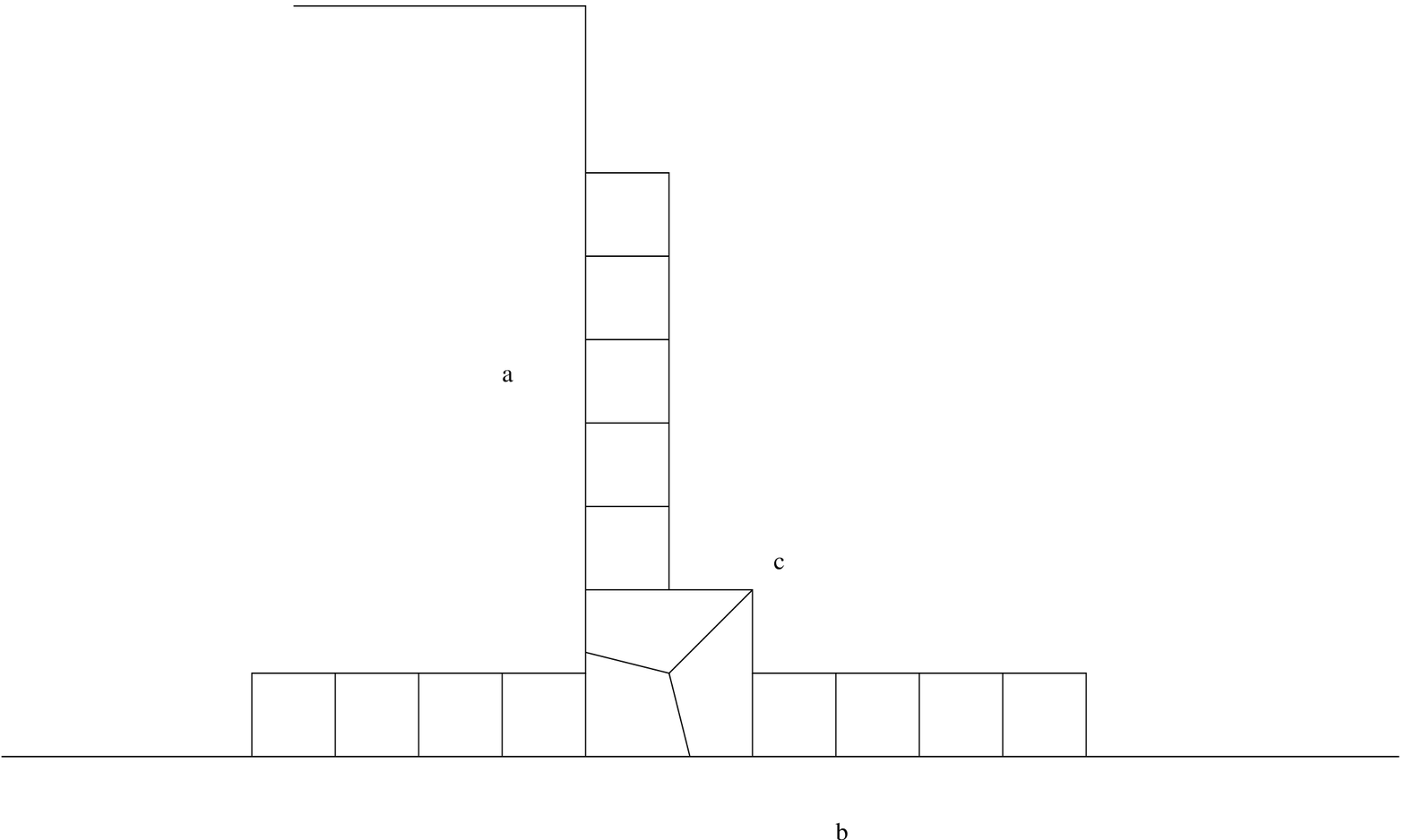, height = 5 cm}
\end{array}$$
\caption{A cartoon of a thin triangle in the decagon tiling.}
\label{Fig:Triangle}
\end{figure}

If $\gamma_3$ contains an end and a side of a corridor $C$, forming a
bad turn, then sweep $\gamma_3$ across $C$.  This occurs at most
twice. Call the resulting curve, which must be well-tempered,
$\gamma$.  These four stages simply move $\alpha \cdot \beta$ through
the thin triangle bounded by $\alpha$, $\beta$, and $\gamma$.  See
Figure~\ref{Fig:Triangle}.

Recall that $\calM$ is the set of words corresponding to paths in the
one-skeleton of the tiling $\calD$.  For any compressed word $\AA$ in
$\calM$ let $\gamma_A$ be the corresponding path.  Again $\gamma_A$ is
only defined up to the action of $\pi_1(S)$.  The above discussion
leads both to a proof of Theorem~\ref{Thm:WellTempered} and of:

\begin{theorem}
\label{Thm:WellTemperedAlgorithm}
There is a polynomial-time algorithm that, given a compressed word
$\AA$ in $\calM$, finds a compressed word $\XX$ where
\begin{itemize}
\item
the path $\gamma_X$ is well-tempered and 
\item
$\gamma_A$ and $\gamma_X$ are homotopic rel endpoints.
\end{itemize}
\end{theorem}

\begin{proof}[Proof sketch]
Suppose that the root of $\AA$ has production $A \produces B \cdot C$.
By induction assume that $\gamma_B$ and $\gamma_C$ are well-tempered.
The first stage is a straight-forward application of Plandowski's
Algorithm (Theorem~\ref{Thm:PlandowskiII}).  The second and fourth
stages require both Plandowski's Algorithm and the fact that words of
period four are highly compressible. To deal with the third simply
examine a constant sized suffix of $w(B'')$ and a constant sized
prefix of $w(C'')$, where $B''$ and $C''$ are the compressed words
output by the second stage.
\end{proof}

\begin{remark}
\label{Rem:Size}
There is a subtlety hidden in this proof sketch -- the compressed word
$\XX$ produced may have bit-size larger than that of $\AA$.  Since the
proof is inductive the growth must be carefully controlled, in part
using Hagenah's Algorithm (Theorem~\ref{Thm:Hagenah})
\end{remark}


\begin{remark}
\label{Rem:ShortLexII}
It may be possible to prove versions of
Theorems~\ref{Thm:WellTempered} and~\ref{Thm:WellTemperedAlgorithm}
using short-lex paths.  However the number of combinatorial
possibilities appears to greatly increase.  Also, I do not know how to
control the growth in size indicated in Remark~\ref{Rem:Size} when
using short-lex paths.  If this could be done then the entire
discussion should apply to general word hyperbolic groups.
\end{remark}

From Theorem~\ref{Thm:WellTemperedAlgorithm} deduce:

\begin{corollary}
\label{Cor:WordProbInSurface}
The word problem for compressed words in $\pi_1(S)$ is solvable in
polynomial time.  This gives a solution to the word problem in
$\Aut(\pi_1(S))$.  \qed
\end{corollary}

The compressed conjugacy problem follows from a careful reading of
Epstein and Holt's paper~\cite{EpsteinHolt05}.  In a few places their
subroutines, acting on words, must be altered to act on compressed
words.  In particular a solution to the fully compressed matching
problem (see the proof of Theorem~\ref{Thm:ConjProbInFree}) replaces
the Knuth-Morris-Pratt algorithm for checking if two words are cyclic
conjugates and for computing roots.  Also, as noted above, the
language of well-tempered paths is regular; this is used in their
proof to find certain bounds.  Thus:

\begin{theorem}
\label{Thm:ConjProbInSurface}
The conjugacy problem for compressed words in $\pi_1(S)$ is solvable
in polynomial time. \qed
\end{theorem}

Since $\pi_1(S)$ has no torsion the simple version of Bridson and
Howie's algorithm~\cite{BridsonHowie05}, adapted to compressed words,
now solves the membership problem for the inner automorphism group
$\Inn(\pi_1(S))$.  Finally, recall Nielsen's Theorem (see \cite[page
175]{MagnusEtAl66}): the mapping class group $\MCG(S)$ is isomorphic
to the outer automorphism group of $\pi_1(S)$.  So, similar to the
proof of \refthm{Out}, well-tempered paths give:

\begin{theorem}[\cite{Mosher95}, \cite{Hamidi-Tehrani00}]
\label{Thm:MCG}
The word problem in $\MCG(S)$ is solvable in polynomial time. \qed
\end{theorem}

\bibliographystyle{plain}
\bibliography{bibfile}
\end{document}

%% file: header.tex
\usepackage{color}
\definecolor{darkgreen}{rgb}{0, 0.40, 0}

\usepackage[
ps2pdf,                
colorlinks=true,
linkcolor=blue,
citecolor=darkgreen,
urlcolor=magenta
]{hyperref}

\usepackage{epsfig}
\usepackage{latexsym}
\usepackage{psfrag}
\usepackage{amssymb}
\usepackage{amsmath}
\usepackage{amsthm} 
\usepackage{verbatim} 

\newcommand{\etal}{{\it et\thinspace al}}

\renewcommand{\th}{{\rm th}}

\newcommand{\calA}{{\mathcal{A}}}

\newcommand{\calD}{{\mathcal{D}}}

\newcommand{\calH}{{\mathcal{H}}}

\newcommand{\calL}{{\mathcal{L}}}
\newcommand{\calM}{{\mathcal{M}}}

\newcommand{\calP}{{\mathcal{P}}}
\newcommand{\calQ}{{\mathcal{Q}}}

\newcommand{\calX}{{\mathcal{X}}}

\renewcommand{\AA}{{\mathbb{A}}}
\newcommand{\BB}{{\mathbb{B}}}

\newcommand{\DD}{{\mathbb{D}}}

\newcommand{\FF}{{\mathbb{F}}}

\newcommand{\NN}{{\mathbb{N}}}

\newcommand{\WW}{{\mathbb{W}}}
\newcommand{\XX}{{\mathbb{X}}}

\renewcommand{\setminus}{{\smallsetminus}}

\newcommand{\Implies}{\rightarrow}

\newcommand{\st}{\mathbin{\mid}} 
\newcommand{\from}{\colon} 

\newcommand{\isom}{{\mathrel{\cong}}} 

\newcommand{\bdy}{{\partial}} 

\newcommand{\MCG}{{\mathcal{MCG}}} 

\newcommand{\Aut}{{\operatorname{Aut}}} 
\newcommand{\Out}{{\operatorname{Out}}} 
\newcommand{\Inn}{{\operatorname{Inn}}} 

\newcommand{\refthm}[1]{Theorem~\ref{Thm:#1}}

\newcommand{\reffig}[1]{Figure~\ref{Fig:#1}}

\theoremstyle{plain}
\newtheorem{theorem}{Theorem}[section]
\newtheorem{corollary}[theorem]{Corollary}
\newtheorem{lemma}[theorem]{Lemma}

\theoremstyle{definition}

\newtheorem{remark}[theorem]{Remark}
\newtheorem{claim}[theorem]{Claim}

\newtheorem{example}[theorem]{Example}

\newcommand{\fakeenv}{} 

\newenvironment{restate}[2]  
{ 
 \renewcommand{\fakeenv}{#2} 
 \theoremstyle{plain} 
 \newtheorem*{\fakeenv}{#1~\ref{#2}} 
 \begin{\fakeenv}
}
{
 \end{\fakeenv}
}





%% file: word.bbl
\begin{thebibliography}{10}

\bibitem{AgolEtAl02}
Ian Agol, Joel Hass, and William Thurston.
\newblock 3-manifold knot genus is {NP}-complete.
\newblock In {\em Proceedings of the Thirty-Fourth Annual ACM Symposium on
  Theory of Computing}, pages 761--766 (electronic), New York, 2002. ACM.
\newblock
  \href{http://front.math.ucdavis.edu/math.GT/0205057}{arXiv:math.GT/0205057}.

\bibitem{Baumslag02}
Gilbert Baumslag, Alexei~G. Myasnikov, and Vladimir Shpilrain.
\newblock Open problems in combinatorial group theory. {S}econd edition.
\newblock In {\em Combinatorial and geometric group theory (New York,
  2000/Hoboken, NJ, 2001)}, volume 296 of {\em Contemp. Math.}, pages 1--38.
  Amer. Math. Soc., Providence, RI, 2002.
\newblock
  \href{http://www.sci.ccny.cuny.edu/~shpil/gworld/problems/oproblems.html}{ht%
tp://www.sci.ccny.cuny.edu/{$\sim$}shpil/gworld/problems/oproblems.html}.

\bibitem{Bigelow01}
Stephen~J. Bigelow.
\newblock Braid groups are linear.
\newblock {\em J. Amer. Math. Soc.}, 14(2):471--486 (electronic), 2001.
\newblock
  \href{http://www.ams.org/joursearch/servlet/DoSearch?f4=author&v4=Bigelow}{h%
ttp://www.ams.org/}.

\bibitem{BirmanEtAl98}
Joan Birman, Ki~Hyoung Ko, and Sang~Jin Lee.
\newblock A new approach to the word and conjugacy problems in the braid
  groups.
\newblock {\em Adv. Math.}, 139(2):322--353, 1998.
\newblock
  \href{http://www.math.columbia.edu/~jb/papers.html}{http://www.math.columbia%
.edu/{$\sim$}jb/papers.html}.

\bibitem{BridsonGroves05}
Martin~R. Bridson and Daniel Groves.
\newblock {Free-group automorphisms, train tracks and the beaded
  decomposition}.
\newblock
  \href{http://front.math.ucdavis.edu/math.GR/0507589}{arXiv:math.GR/0507589}.

\bibitem{BridsonHowie05}
Martin~R. Bridson and James Howie.
\newblock {Conjugacy of finite subsets in hyperbolic groups}.
\newblock
  \href{http://www.ma.ic.ac.uk/~mbrids/papers/bhowie/}{http://www.ma.ic.ac.uk/%
{$\sim$}mbrids/papers/bhowie/}.

\bibitem{EpsteinHolt05}
David Epstein and Derek Holt.
\newblock The linearity of the conjugacy problem in word-hyperbolic groups.
\newblock 2005.
\newblock
  \href{http://www.maths.warwick.ac.uk/~dfh/download/papers/}{http://www.maths%
.warwick.ac.uk/{$\sim$}dfh/download/papers/}.

\bibitem{EpsteinEtAl92}
David B.~A. Epstein, James~W. Cannon, Derek~F. Holt, Silvio V.~F. Levy,
  Michael~S. Paterson, and William~P. Thurston.
\newblock {\em Word processing in groups}.
\newblock Jones and Bartlett Publishers, Boston, MA, 1992.

\bibitem{GKPR96}
Leszek Gasieniec, Marek Karpinski, Wojciech Plandowski, and Wojciech Rytter.
\newblock Efficient algorithms for {L}empel-{Z}iv encoding.
\newblock 1097:392--403, 1996.
\newblock
  \href{http://citeseer.ist.psu.edu/22169.html}{http://citeseer.ist.psu.edu/22%
169.html}.

\bibitem{Hagenah00}
Christian Hagenah.
\newblock {\em {\em Gleichungen mit regul{\"a}ren Randbedingungen {\"u}ber
  freien Gruppen}}.
\newblock Dissertation, Universit{\"a}t Stuttgart, Fakult{\"a}t Informatik,
  Elektrotechnik und Informationstechnik, August 2000.
\newblock \href{http://elib.uni-stuttgart.de/opus/volltexte/2000/673/}
  {http://elib.uni-stuttgart.de/opus/volltexte/2000/673/}.

\bibitem{Hamidi-Tehrani00}
Hessam Hamidi-Tehrani.
\newblock On complexity of the word problem in braid groups and mapping class
  groups.
\newblock {\em Topology Appl.}, 105(3):237--259, 2000.
\newblock
  \href{http://www.math.columbia.edu/~hessam/}{http://www.math.columbia.edu/{$%
\sim$}hessam/}.

\bibitem{Hempel76}
John Hempel.
\newblock {\em $3$-{M}anifolds}.
\newblock Princeton University Press, Princeton, N. J., 1976.
\newblock Ann. of Math. Studies, No. 86.

\bibitem{Humphries79}
Stephen~P. Humphries.
\newblock Generators for the mapping class group.
\newblock In {\em Topology of low-dimensional manifolds (Proc. Second Sussex
  Conf., Chelwood Gate, 1977)}, volume 722 of {\em Lecture Notes in Math.},
  pages 44--47. Springer, Berlin, 1979.

\bibitem{KarpinskiEtAl95}
Marek Karpinski, Wojciech Rytter, and Ayumi Shinohara.
\newblock Pattern-matching for strings with short descriptions.
\newblock In {\em Combinatorial pattern matching (Espoo, 1995)}, volume 937 of
  {\em Lecture Notes in Comput. Sci.}, pages 205--214. Springer, Berlin, 1995.

\bibitem{Krammer02}
Daan Krammer.
\newblock Braid groups are linear.
\newblock {\em Ann. of Math. (2)}, 155(1):131--156, 2002.
\newblock
  \href{http://www.maths.warwick.ac.uk/~daan/}{http://www.maths.warwick.ac.uk/%
{$\sim$}daan/}.

\bibitem{Lohrey04}
Markus Lohrey.
\newblock Word problems on compressed words.
\newblock In {\em Automata, languages and programming}, volume 3142 of {\em
  Lecture Notes in Comput. Sci.}, pages 906--918. Springer, Berlin, 2004.
\newblock
  \href{http://www.informatik.uni-stuttgart.de/fmi/ti/personen/Lohrey/veroeff_%
eng.html}{http://www.informatik.uni-stuttgart.de/fmi/ti/personen/Lohrey/}.

\bibitem{Lothaire83}
M.~Lothaire.
\newblock {\em Combinatorics on words}, volume~17 of {\em Encyclopedia of
  Mathematics and its Applications}.
\newblock Addison-Wesley Publishing Co., Reading, Mass., 1983.
\newblock A collective work by Dominique Perrin, Jean Berstel, Christian
  Choffrut, Robert Cori, Dominique Foata, Jean Eric Pin, Guiseppe Pirillo,
  Christophe Reutenauer, Marcel-P. Sch\"utzenberger, Jacques Sakarovitch and
  Imre Simon, With a foreword by Roger Lyndon, Edited and with a preface by
  Perrin.

\bibitem{LyndonSchupp77}
Roger~C. Lyndon and Paul~E. Schupp.
\newblock {\em Combinatorial group theory}.
\newblock Springer-Verlag, Berlin, 1977.
\newblock Ergebnisse der Mathematik und ihrer Grenzgebiete, Band 89.

\bibitem{MagnusEtAl66}
Wilhelm Magnus, Abraham Karrass, and Donald Solitar.
\newblock {\em Combinatorial group theory: {P}resentations of groups in terms
  of generators and relations}.
\newblock Interscience Publishers [John Wiley \& Sons, Inc.], New
  York-London-Sydney, 1966.

\bibitem{Miyazaki00}
Masamichi Miyazaki, Ayumi Shinohara, and Masayuki Takeda.
\newblock An improved pattern matching algorithm for strings in terms of
  straight-line programs.
\newblock {\em J. Discrete Algorithms (Oxf.)}, 1(1):187--204, 2000.
\newblock
  \href{http://www.shino.ecei.tohoku.ac.jp/~ayumi/publications.html}{http://ww%
w.shino.ecei.tohoku.ac.jp/{$\sim$}ayumi/publications.html}.

\bibitem{Mosher95}
Lee Mosher.
\newblock Mapping class groups are automatic.
\newblock {\em Ann. of Math. (2)}, 142(2):303--384, 1995.

\bibitem{Plandowski94}
Wojciech Plandowski.
\newblock Testing equivalence of morphisms on context-free languages.
\newblock In {\em Algorithms---ESA '94 (Utrecht)}, volume 855 of {\em Lecture
  Notes in Comput. Sci.}, pages 460--470. Springer, Berlin, 1994.

\bibitem{SchaeferEtAl02}
Marcus Schaefer, Eric Sedgwick, and Daniel {\v{S}}tefankovi{\v{c}}.
\newblock Algorithms for normal curves and surfaces.
\newblock In {\em Computing and combinatorics}, volume 2387 of {\em Lecture
  Notes in Comput. Sci.}, pages 370--380. Springer, Berlin, 2002.
\newblock
  \href{http://www.cs.rochester.edu/~stefanko/}{http://www.cs.rochester.edu/{$%
\sim$}stefanko/}.

\end{thebibliography}
